\documentclass[11pt,reqno]{amsart}
\usepackage[T1]{fontenc}
\usepackage{lmodern}
\usepackage{amsmath,amssymb,amsthm,mathtools}
\usepackage{microtype}
\usepackage[hypertexnames=false,colorlinks=true,linkcolor=black,citecolor=black,urlcolor=blue]{hyperref}
\allowdisplaybreaks[2]

\newcommand{\C}{\mathbb{C}}
\newcommand{\Z}{\mathbb{Z}}

\newcommand{\dq}{D_q}
\newcommand{\Xq}{\mathbb{R}_{q,+}}
\newcommand{\Sq}{\mathcal{S}_q}
\newcommand{\ind}{\mathbf{1}}
\newcommand{\e}{\mathrm{e}}
\newcommand{\qbinom}[2]{\genfrac{[}{]}{0pt}{}{#1}{#2}_q}

\newtheorem{theorem}{Theorem}[section]
\newtheorem{proposition}[theorem]{Proposition}
\theoremstyle{definition}
\newtheorem{definition}[theorem]{Definition}
\theoremstyle{remark}
\newtheorem{remark}[theorem]{Remark}
\numberwithin{equation}{section}
\hypersetup{pdftitle={q-Completely monotonic and q-Bernstein functions},pdfauthor={Valmir Krasniqi}}
\title{q-Completely monotonic and q-Bernstein functions}
\author{Valmir Krasniqi}
\address{Institute of Science, Technology, Engineering and Mathematics (STEM), Rr. Nekibe Kelmendi Nr. 26, 10000 Prishtin\"e, Republic of Kosovo}
\email{valmir@stem-ks.org}
\date{}
\begin{document}
\begin{abstract}
For $0<q<1$, we study $q$-completely monotone and $q$-Bernstein functions defined by the Jackson derivative.
We establish the positive linear closure of the latter class and prove an integral characterization of $q$-Bernstein functions continuous at zero on a geometric lattice.
The representing measures satisfy the natural L\'evy integrability condition, and the corresponding integral kernels are positive semidefinite.
We prove that classical completely monotone and Bernstein functions belong to the respective $q$-classes and examine composition with absolutely monotone functions.
Explicit examples show the limits of product closure, composition, and definiteness with respect to ordinary addition.
\end{abstract}
\maketitle

\section{Preliminaries}

Throughout the paper, $0<q<1$ is fixed.
We use the standard notation of $q$-calculus; see \cite{Chung,Ernst2000,ErnstBook,Exton}.
For $a\in\C$ and integers $n\geq0$, $0\leq k\leq n$, put
\[
 [a]_q=\dfrac{1-q^a}{1-q},\qquad
 [n]_q!=\prod_{j=1}^{n}[j]_q,\qquad [0]_q!=1,
 \qquad
 \qbinom{n}{k}=\dfrac{[n]_q!}{[k]_q![n-k]_q!}.
\]
The two $q$-exponential functions, in the normalization used here, are
\begin{equation}\label{eq:exponentials}
 e_q(z)=\sum_{n=0}^{\infty}\dfrac{z^n}{[n]_q!},
 \qquad
 E_q(z)=\sum_{n=0}^{\infty}q^{n(n-1)/2}\dfrac{z^n}{[n]_q!}.
\end{equation}
The first series has radius of convergence $(1-q)^{-1}$, whereas the second is entire.
With $(z;q)_\infty=\prod_{j=0}^{\infty}(1-zq^j)$, their product representations are
\begin{equation}\label{eq:products}
 e_q(z)=\dfrac{1}{((1-q)z;q)_\infty},
 \qquad E_q(z)=(-(1-q)z;q)_\infty.
\end{equation}
The product on the right extends $e_q$ to the negative real axis.

A set $I\subseteq(0,\infty)$ is called $q$-stable if $qx\in I$ whenever $x\in I$.
For a finite-valued function on such a set, the Jackson derivative is
\begin{equation}\label{eq:jackson}
 \dq f(x)=\dfrac{f(x)-f(qx)}{(1-q)x},\qquad x\in I.
\end{equation}
For each fixed $x$ and $n$, the iterate $\dq^n f(x)$ is a finite linear combination of the values $f(x),f(qx),\ldots,f(q^nx)$.
In particular, the existence of these iterates does not require ordinary differentiability.
The standard differentiation identities are
\begin{equation}\label{eq:expderivatives}
 \dq^n e_q(ax)=a^n e_q(ax),\qquad
 \dq^n E_q(ax)=a^n q^{n(n-1)/2}E_q(aq^nx).
\end{equation}
The first identity holds in the Taylor disk and wherever its product continuation and the relevant difference quotients are defined.
The second holds for every complex argument.
We shall also use the Jackson integral
\begin{equation}\label{eq:jacksonintegral}
 \int_0^x h(y)\,\mathrm d_qy
 =(1-q)x\sum_{j=0}^{\infty}q^j h(q^jx),
\end{equation}
whenever the series converges.

A nonnegative function $g:I\to[0,\infty)$ is \emph{$q$-completely monotone}, or $q$-CM, if
\begin{equation}\label{eq:qcm}
 (-1)^n\dq^n g(x)\geq0\qquad(x\in I,\ n=0,1,2,\ldots).
\end{equation}
The $q$-discrete theory of complete monotonicity and its $q$-Laplace transform are developed in \cite[Section~5]{Dhaouadi}.
The standard $q$-Leibniz formula
\[
 \dq^n(gh)(x)=\sum_{k=0}^n\qbinom{n}{k}
 (\dq^{n-k}g)(q^kx)(\dq^kh)(x)
\]
shows that the class of $q$-CM functions is closed under nonnegative linear combinations and products.
These recalled facts will be used without separate proofs.

For the integral representations involving $E_q$, we distinguish the geometric domain and the spectral lattice:
\begin{equation}\label{eq:lattices}
 \Xq=\{q^m:m\in\Z\},\qquad
 \Sq=\left\{\dfrac{q^k}{1-q}:k\in\Z\right\}.
\end{equation}
Both sets are invariant under multiplication by $q$ and by $q^{-1}$.
At $0$, a function on $\Xq$ is assigned its finite lattice limit when that limit exists.
For $t\in\Sq$, write
\[
 K_t(x)=E_q(-tx),\qquad x\in\Xq,\qquad K_t(0)=1.
\]
All spectral integrals below are ordinary integrals against positive Borel measures on the indicated lattice; they are not additional Jackson integrals.
No boundedness assumption is imposed on a density of these measures with respect to Jackson measure.

For comparison, a classical completely monotone function $g$ is a nonnegative $C^\infty$ function on $(0,\infty)$ satisfying $(-1)^ng^{(n)}\geq0$ for every $n\geq0$.
A classical Bernstein function is a nonnegative $C^\infty$ function $f$ with completely monotone derivative.
We recall, without proof, its L\'evy--Khintchine representation \cite[Theorem~3.2]{Schilling}:
\begin{equation}\label{eq:classicallevy}
 f(x)=b+ax+\int_{(0,\infty)}(1-\e^{-tx})\,\mathrm d\mu(t),
 \qquad
 \int_{(0,\infty)}\min(1,t)\,\mathrm d\mu(t)<\infty,
\end{equation}
where $a,b\geq0$ and $\mu$ is positive.
Background on these classical classes and their composition properties can be found in \cite{MillerSamko,Schilling,Srivastava}.

A kernel $K$ on a set $X$ is called \emph{positive definite} (in the positive semidefinite sense) if
\[
 \sum_{i,j=1}^{m}c_i\overline{c_j}K(x_i,x_j)\geq0
 \quad(m\geq1,\ x_i\in X,\ c_i\in\C).
\]
For a real function $u$ on $[0,\infty)$, positive definiteness with respect to ordinary addition means that $K(x,y)=u(x+y)$ has this property.
Negative definiteness in the semigroup sense used here means that
\begin{equation}\label{eq:negativekernel}
 N_u(x,y)=u(x)+u(y)-u(x+y)
\end{equation}
is positive definite; see \cite[Chapter~4]{Schilling}.
These notions use ordinary addition, which is distinct from the formal $q$-addition considered below.

\section{Main results}

\begin{definition}\label{def:qbf}
A nonnegative function $f$ on a $q$-stable set $I\subseteq(0,\infty)$ is called a \emph{$q$-Bernstein function} if
\[
 (-1)^{n-1}\dq^n f(x)\geq0\qquad(x\in I,\ n=1,2,\ldots).
\]
Equivalently, $\dq f$ is $q$-CM on $I$.
\end{definition}

\begin{theorem}\label{th:closure}
Let $f_j$ ($1\leq j\leq m$) be $q$-Bernstein functions on $I$.
If $a_j\geq0$ for $1\leq j\leq m$, then $\sum_{j=1}^m a_jf_j$ is $q$-Bernstein.
This class is not closed under pointwise multiplication.
\end{theorem}
\begin{proof}
For $n\geq1$, linearity gives
\[
 (-1)^{n-1}\dq^n\left(\sum_{j=1}^m a_jf_j\right)
 =\sum_{j=1}^m a_j(-1)^{n-1}\dq^nf_j\geq0.
\]
Nonnegativity of the sum is immediate.
For the product assertion, both functions $f_1(x)=f_2(x)=1+x$ are $q$-Bernstein, since $\dq f_j=1$ and $\dq^nf_j=0$ for $n\geq2$.
However,
\[
 \dq^2\bigl((1+x)^2\bigr)=1+q>0,
\]
which violates the required second-order inequality.
\end{proof}

\begin{theorem}[Spectral representation]\label{th:representation}
Let $f:\Xq\cup\{0\}\to[0,\infty)$ be continuous at $0$ along the lattice.
The following three conditions are equivalent.
\begin{enumerate}
\item[\textup{(i)}] There exist $a,b\geq0$ and a positive measure $\mu$ on $\Sq$ such that
\begin{equation}\label{eq:qlevy}
 f(x)=b+ax+\int_{\Sq}(1-K_t(x))\,\mathrm d\mu(t),
 \qquad
 \int_{\Sq}\min(1,t)\,\mathrm d\mu(t)<\infty.
\end{equation}
\item[\textup{(ii)}] There exist $a\geq0$ and a positive measure $\nu$ on $\Sq$ such that
\begin{equation}\label{eq:derivativerep}
 \dq f(x)=a+\int_{\Sq}K_t(x)\,\mathrm d\nu(t),
 \qquad
 \int_{\Sq}\min(1,t^{-1})\,\mathrm d\nu(t)<\infty.
\end{equation}
\item[\textup{(iii)}] The function $f$ is $q$-Bernstein on $\Xq$.
\end{enumerate}
In these representations, $b=f(0)$.
For the corresponding representations, the measures are related by
\begin{equation}\label{eq:measuremaps}
 \nu(A)=\int_{\Sq}t\ind_A(qt)\,\mathrm d\mu(t),
 \qquad
 \mu(A)=\int_{\Sq}\dfrac{q}{t}\ind_A(t/q)\,\mathrm d\nu(t)
\end{equation}
for every Borel set $A\subseteq\Sq$.
\end{theorem}

\begin{proof}
We first record the positivity and bounds needed for the implications below.
If $t=q^k/(1-q)$ and $x=q^m$, then
\[
 K_t(x)=\prod_{j=0}^{\infty}(1-q^{k+m+j}).
\]
When $k+m\leq0$, one factor is zero, so $K_t(x)=0$.
When $k+m\geq1$, every factor belongs to $(0,1)$.
Consequently $0\leq K_t(x)\leq1$ on the two lattices.
Writing $u=(1-q)tx$, the product also gives
\begin{equation}\label{eq:kernelbounds}
 \min(1,(1-q)tx)\leq1-K_t(x)\leq\min(1,tx).
\end{equation}
Indeed, when $u\leq1$, the lower bound follows from $K_t(x)\leq1-u$, while the upper bound follows from
\[
 1-\prod_{j=0}^{\infty}(1-uq^j)
 \leq\sum_{j=0}^{\infty}uq^j=tx.
\]
When $u>1$, the lattice condition implies $K_t(x)=0$, and both bounds follow directly.
In addition, for every $n\geq0$,
\begin{equation}\label{eq:kernelderivatives}
 (-1)^n\dq^nK_t(x)
 =t^n q^{n(n-1)/2}K_t(q^nx)\geq0.
\end{equation}
For fixed $x$ and $n$, the kernel on the right vanishes whenever
\begin{equation}\label{eq:spectralcutoff}
 t\geq\dfrac{q^{-n}}{(1-q)x}.
\end{equation}

Suppose first that \textup{(i)} holds.
The upper bound in \eqref{eq:kernelbounds} makes the integral in \eqref{eq:qlevy} finite.
For $0<x\leq x_0$ along the lattice, its integrand is bounded by the integrable function $\min(1,tx_0)$ and tends to zero as $x\to0$.
Dominated convergence therefore gives $f(0)=b$.
Each Jackson iterate is a finite linear combination of evaluations, so it may be applied under the integral in \eqref{eq:qlevy}.
Using \eqref{eq:kernelderivatives}, we obtain
\begin{equation}\label{eq:levyderivatives}
 (-1)^{n-1}\dq^nf(x)
 =a\ind_{\{1\}}(n)
 +q^{n(n-1)/2}\int_{\Sq}t^nK_t(q^nx)\,\mathrm d\mu(t),
 \qquad n\geq1.
\end{equation}
These integrals are finite: by \eqref{eq:spectralcutoff} the integrands vanish above a finite threshold, and on any bounded interval $t^n$ is bounded by a constant times $\min(1,t)$ for $n\geq1$.
Every term in \eqref{eq:levyderivatives} is nonnegative, proving the $q$-Bernstein inequalities.
For $n=1$, the formula becomes
\[
 \dq f(x)=a+\int_{\Sq}tK_{qt}(x)\,\mathrm d\mu(t).
\]
The first transformation in \eqref{eq:measuremaps} now gives \eqref{eq:derivativerep}, and its integrability condition follows from
\[
 \int_{\Sq}\min(1,s^{-1})\,\mathrm d\nu(s)
 =\int_{\Sq}\min(t,q^{-1})\,\mathrm d\mu(t)<\infty.
\]

Conversely, assume \textup{(ii)}.
For fixed $x$, the integral in \eqref{eq:derivativerep} is finite because $K_t(x)$ vanishes above a finite threshold and $\nu((0,T]\cap\Sq)<\infty$ for every $T<\infty$.
The last assertion follows directly from the weighted integrability in \textup{(ii)}.
The finite telescoping identity
\[
 (1-q)x\sum_{j=0}^{N}q^j\dq f(q^jx)
 =f(x)-f(q^{N+1}x)
\]
and continuity at zero yield
\begin{align}\label{eq:integratedderivative}
 f(x)-f(0)
 &=\int_0^x\dq f(y)\,\mathrm d_qy\notag\\
 &=ax+\int_{\Sq}\left(\int_0^xK_t(y)\,\mathrm d_qy\right)\mathrm d\nu(t).
\end{align}
Tonelli's theorem applies because the kernels on both lattices are nonnegative.
The exact Jackson antiderivative is
\begin{equation}\label{eq:antiderivative}
 \int_0^xK_t(y)\,\mathrm d_qy
 =\dfrac{q}{t}\left(1-K_{t/q}(x)\right).
\end{equation}
To check this identity, \eqref{eq:expderivatives} shows that the Jackson derivative of the right-hand side is $K_t(x)$; its limit at zero is zero, so the telescoping identity gives \eqref{eq:antiderivative}.
Define $\mu$ by the second transformation in \eqref{eq:measuremaps}.
Because $t/q\in\Sq$, substitution in \eqref{eq:integratedderivative} gives \eqref{eq:qlevy} with $b=f(0)$.
Moreover,
\[
 \int_{\Sq}\min(1,s)\,\mathrm d\mu(s)
 =\int_{\Sq}\min(1,q/t)\,\mathrm d\nu(t)<\infty.
\]
The two transformations in \eqref{eq:measuremaps} are inverse: applying both successively multiplies an integrand by $t\cdot q/(qt)=1$ and returns the spectral variable to its original value.
Thus \textup{(i)} and \textup{(ii)} are equivalent and imply \textup{(iii)}.

It remains to derive \textup{(ii)} from \textup{(iii)}.
We first reconstruct the spectral measure of an arbitrary $q$-CM function $g$ on $\Xq$.
The inequality $\dq g\leq0$ implies that $g(q^m)$ is nondecreasing in $m$.
Hence
\[
 a=\lim_{m\to-\infty}g(q^m)\geq0,
 \qquad h_m=g(q^m)-a
\]
are well defined, $h_m\geq0$, and $h_m\to0$ as $m\to-\infty$.
Let $\mathsf T$ denote the backward shift, $(\mathsf T h)_m=h_{m-1}$, and let $\mathsf I$ be the identity operator.
Iteration of \eqref{eq:jackson} gives, for $n\geq1$,
\begin{equation}\label{eq:positivefiniteproduct}
 \left[\prod_{r=0}^{n-1}(\mathsf I-q^r\mathsf T)h\right]_m
 =(-1)^n(1-q)^nq^{nm-n(n+1)/2}
   \dq^ng(q^{m-n})\geq0.
\end{equation}
Put $(q;q)_j=\prod_{r=1}^j(1-q^r)$, with $(q;q)_0=1$, and $C_q=(q;q)_\infty>0$.
The finite $q$-binomial formula and passage to the limit in \eqref{eq:positivefiniteproduct} yield
\begin{equation}\label{eq:spectralcoefficients}
 v_m=\sum_{j=0}^{\infty}
 \dfrac{(-1)^jq^{j(j-1)/2}}{(q;q)_j}\,h_{m-j}\geq0.
\end{equation}
Indeed, the coefficient at order $j$ in the finite product is
$(-1)^jq^{j(j-1)/2}\qbinom{n}{j}$, which converges to the coefficient in \eqref{eq:spectralcoefficients} and has absolute value at most $q^{j(j-1)/2}/(q;q)_j$.
Since $0\leq h_{m-j}\leq h_m$ and the absolute coefficients are summable, dominated convergence justifies the limit.

Write $p_m=h_m-h_{m-1}$.
For each fixed $M$, both $h$ and $p$ are bounded on the backward half-lattice $m\leq M$.
The standard Euler product and reciprocal identities \cite{ErnstBook,Exton} therefore apply as absolutely convergent series in the backward shift:
\[
 v=(q\mathsf T;q)_\infty p,
 \qquad
 (q\mathsf T;q)_\infty^{-1}
 =\sum_{j=0}^{\infty}\dfrac{q^j}{(q;q)_j}\mathsf T^j.
\]
Here the operator products are defined by their power-series coefficients.
Those of both $(qz;q)_\infty$ and its reciprocal are absolutely summable for $|z|\leq1$, so multiplication of the series is justified on each bounded half-lattice.
It follows that
\[
 p_m=\sum_{j=0}^{\infty}\dfrac{q^j}{(q;q)_j}v_{m-j}.
\]
Because $h_m=\sum_{r=0}^{\infty}p_{m-r}$ and $v_m\geq0$, Tonelli's theorem gives
\begin{align}\label{eq:latticereconstruction}
 h_m
 &=\sum_{\ell=0}^{\infty}v_{m-\ell}
   \sum_{j=0}^{\ell}\dfrac{q^j}{(q;q)_j}
 =\sum_{\ell=0}^{\infty}\dfrac{v_{m-\ell}}{(q;q)_\ell}.
\end{align}
The last identity uses
\[
 \sum_{j=0}^{\ell}\dfrac{q^j}{(q;q)_j}
 =\dfrac{1}{(q;q)_\ell},
\]
which follows by telescoping.
Define a positive measure on $\Sq$ by
\begin{equation}\label{eq:reconstructedmeasure}
 \nu\left(\left\{\dfrac{q^k}{1-q}\right\}\right)
 =\dfrac{v_{1-k}}{C_q},\qquad k\in\Z.
\end{equation}
For $k+m\geq1$, one has
\[
 K_{q^k/(1-q)}(q^m)=\dfrac{C_q}{(q;q)_{k+m-1}},
\]
and this kernel is zero for $k+m\leq0$.
Equation~\eqref{eq:latticereconstruction} consequently becomes
\begin{equation}\label{eq:generalcmrepresentation}
 g(x)=a+\int_{\Sq}K_t(x)\,\mathrm d\nu(t),\qquad x\in\Xq.
\end{equation}
These integrals are finite by the reconstruction itself.
Moreover, $K_t(q^m)\geq C_q$ whenever $t\leq q^{1-m}/(1-q)$, so finiteness of $g(q^m)$ implies
$\nu((0,T]\cap\Sq)<\infty$ for every finite $T$.

Now assume \textup{(iii)} and apply this reconstruction to $g=\dq f$, which is $q$-CM by definition.
The telescoping and Tonelli calculation in \eqref{eq:integratedderivative} uses only positivity, finite pointwise integrals and continuity of $f$ at zero.
It therefore applies to the measure just constructed and gives
\[
 f(x)-f(0)-ax
 =\int_{\Sq}\dfrac{q}{t}(1-K_{t/q}(x))\,\mathrm d\nu(t)<\infty.
\]
For any fixed $x\in\Xq$, the lower bound in \eqref{eq:kernelbounds} implies
\[
 \dfrac{q}{t}(1-K_{t/q}(x))
 \geq\min(q/t,(1-q)x)
 \geq\min(q,(1-q)x)\min(1,t^{-1}).
\]
The positive constant $\min(q,(1-q)x)$ is independent of $t$.
Thus the weighted integrability required in \textup{(ii)} follows from finiteness of the preceding integral.
This proves \textup{(iii)}$\Rightarrow$\textup{(ii)} and completes the proof.
\end{proof}

The corresponding spectral measures in Theorem~\ref{th:representation} satisfy
\[
 \nu(\Sq)=\int_{\Sq}t\,\mathrm d\mu(t),
\]
with either side allowed to be infinite.

\begin{proposition}\label{prop:smallparameter}
Let $f:I\to[0,\infty)$, where $I$ is $q$-stable.
If there exists $\varepsilon>0$ such that $E_q(-kf)$ is $q$-CM on $I$ for every $0<k<\varepsilon$, then $f$ is $q$-Bernstein on $I$.
\end{proposition}
\begin{proof}
Fix $x\in I$ and $n\geq1$.
At each of the finitely many points $q^jx$, $0\leq j\leq n$, the Taylor expansion at zero gives
\[
 E_q(-kf(q^jx))=1-kf(q^jx)+O_{x,n,f}(k^2)
 \qquad(k\downarrow0).
\]
Because $\dq^n$ is a finite linear combination of these values and $\dq^n1=0$, it follows that
\[
 0\leq(-1)^n\dq^n\bigl(E_q(-kf)\bigr)(x)
 =k(-1)^{n-1}\dq^nf(x)+O_{x,n,f}(k^2).
\]
Division by $k>0$ and passage to the limit prove the assertion for every $n\geq1$ and $x\in I$.
\end{proof}

\begin{proposition}\label{prop:positivekernel}
Let $g$ be $q$-CM on $\Xq$.
Then there exist $c\geq0$ and a positive measure $\rho$ on $\Sq$, finite on $(0,T]\cap\Sq$ for every finite $T$, such that
\begin{equation}\label{eq:finitemixture}
 g(x)=c+\int_{\Sq}K_t(x)\,\mathrm d\rho(t),\qquad x\in\Xq.
\end{equation}
The associated kernel
\begin{equation}\label{eq:positivekernel}
 \mathcal K_g(x,y)=c+\int_{\Sq}K_t(x)K_t(y)\,\mathrm d\rho(t)
\end{equation}
is positive definite on $\Xq$.
If $g$ has a finite limit at zero, then $\rho$ is finite and $g(0)=c+\rho(\Sq)$.
\end{proposition}
\begin{proof}
The reconstruction \eqref{eq:spectralcoefficients}--\eqref{eq:generalcmrepresentation} in the proof of Theorem~\ref{th:representation}, applied to $g$, gives \eqref{eq:finitemixture} with $c=\lim_{x\to\infty,\,x\in\Xq}g(x)$.
Since $0\leq K_t\leq1$, the product integral in \eqref{eq:positivekernel} is bounded by the finite integral of $K_t(x)$.
For $x_1,\ldots,x_m\in\Xq$ and $c_1,\ldots,c_m\in\C$,
\[
 \sum_{i,j=1}^m c_i\overline{c_j}\mathcal K_g(x_i,x_j)
 =c\left|\sum_{i=1}^m c_i\right|^2
 +\int_{\Sq}\left|\sum_{i=1}^m c_iK_t(x_i)\right|^2\,\mathrm d\rho(t)
 \geq0.
\]
The integral is finite because $K_t(x_i)^2\leq K_t(x_i)$ and each $K_t(x_i)$ is $\rho$-integrable.
This proves positive definiteness.
Finally, $K_t(x)$ increases to $1$ as $x$ tends to zero along $\Xq$.
Monotone convergence in \eqref{eq:finitemixture} therefore gives $g(0)=c+\rho(\Sq)$ whenever the limit is finite.
\end{proof}

\begin{remark}\label{rem:addition}
The multiplication law for $E_q$ uses formal $q$-addition, not ordinary addition.
For $r>0$, $r\ne1$, define the formal powers by
\[
 (u\oplus_r v)^n
 =\sum_{k=0}^n\genfrac{[}{]}{0pt}{}{n}{k}_r u^{n-k}v^k.
\]
Since $[n]_{1/q}!=q^{-n(n-1)/2}[n]_q!$, one has $E_q=e_{1/q}$ and
\begin{equation}\label{eq:formaladdition}
 E_q(u\oplus_{1/q}v)=E_q(u)E_q(v).
\end{equation}
Here the left side means substitution of the formal powers into the series for $E_q$.
The identity follows by multiplying the two absolutely convergent series and using
\[
 \genfrac{[}{]}{0pt}{}{n}{k}_{1/q}
 =q^{-k(n-k)}\qbinom{n}{k}.
\]
Thus $\mathcal K_g$ is a kernel associated with the specified spectral representation and formal translation.
No scalar binary operation on $\Xq$, or identification with $g(x+y)$, is being assumed.
In particular, \eqref{eq:formaladdition} does not give $E_q(u+v)=E_q(u)E_q(v)$.
\end{remark}

\begin{proposition}\label{prop:monotonicity}
If $f$ is nondecreasing, respectively nonincreasing, on $(0,\infty)$, then $\dq f\geq0$, respectively $\dq f\leq0$.
Strict monotonicity gives the corresponding strict inequality.
\end{proposition}
\begin{proof}
Since $qx<x$ and $(1-q)x>0$, the assertion follows from the sign of the numerator in \eqref{eq:jackson}.
\end{proof}

\begin{theorem}\label{th:classicalcm}
Every classical completely monotone function on $(0,\infty)$ is $q$-CM.
\end{theorem}
\begin{proof}
For any $u\in C^1(0,\infty)$, the fundamental theorem of calculus gives
\begin{equation}\label{eq:oneaverage}
 \dq u(x)=\dfrac{1}{1-q}\int_q^1u'(xs)\,\mathrm ds.
\end{equation}
More generally, for $u\in C^n(0,\infty)$,
\begin{equation}\label{eq:multipleaverage}
 \dq^n u(x)=\dfrac{1}{(1-q)^n}
 \int_{[q,1]^n}
 \left(\prod_{j=1}^n s_j^{\,n-j}\right)
 u^{(n)}(xs_1\cdots s_n)\,\mathrm ds_1\cdots\mathrm ds_n.
\end{equation}
We verify this formula by induction.
The case $n=1$ is \eqref{eq:oneaverage}.
If the formula holds at order $n$, apply $\dq$ under its integral and use \eqref{eq:oneaverage} on the function $x\mapsto u^{(n)}(xP)$, where $P=s_1\cdots s_n$.
This introduces the additional variable $s_{n+1}$ and the factor $P$, changing the weight to
\[
 \left(\prod_{j=1}^n s_j^{n-j}\right)P
 =\prod_{j=1}^{n+1}s_j^{n+1-j}.
\]
All integrals are over compact sets contained in $(0,\infty)$ after multiplication by $x$, so the indicated operations are justified for smooth $u$.
This proves \eqref{eq:multipleaverage} at the next order.
Its weights are nonnegative.
Consequently, if $(-1)^nu^{(n)}\geq0$, then $(-1)^n\dq^nu\geq0$.
Together with $u\geq0$, this is $q$-complete monotonicity.
\end{proof}

\begin{theorem}\label{th:classicalbf}
Every classical Bernstein function on $(0,\infty)$ is a $q$-Bernstein function.
\end{theorem}
\begin{proof}
Such a function $f$ is nonnegative and satisfies $(-1)^{n-1}f^{(n)}\geq0$ for every $n\geq1$.
Applying \eqref{eq:multipleaverage} at each order yields $(-1)^{n-1}\dq^nf\geq0$.
\end{proof}

\begin{remark}\label{rem:counterexample}
The converses of Theorems~\ref{th:classicalcm} and \ref{th:classicalbf} fail in general, even for strictly positive functions smooth at zero.
For an explicit example, set $q=1/10$ and define
\begin{equation}\label{eq:counterexamples}
 B(x)=(1-q)x\sum_{j=0}^{\infty}q^j\e^{-q^jx},
 \qquad F(x)=1-B(x),\qquad H(x)=1+B(x).
\end{equation}
All differentiated series converge uniformly on compact intervals, so these functions belong to $C^\infty[0,\infty)$.
Telescoping gives $\dq B(x)=\e^{-x}$.
Furthermore,
\[
 (1-q)q^jx\e^{-q^jx}
 \leq\int_{q^{j+1}x}^{q^jx}\e^{-t}\,\mathrm dt,
\]
and summation yields $0\leq B(x)\leq1-\e^{-x}<1$ for finite $x\geq0$.
Thus $F(x)>0$ and $H(x)>0$.
By Theorem~\ref{th:classicalcm}, the function $\e^{-x}$ is $q$-CM; hence $F$ is $q$-CM and $H$ is $q$-Bernstein.
However,
\[
 B''(x)=(1-q)\sum_{j=0}^{\infty}q^{2j}\e^{-q^jx}(q^jx-2),
\]
so
\begin{equation}\label{eq:secondderivativebound}
 B''(3)\geq\dfrac{9}{10}\left(\e^{-3}-2\sum_{j=1}^{\infty}q^{2j}\right)
 =\dfrac{9}{10}\left(\e^{-3}-\dfrac{2}{99}\right)>0.
\end{equation}
The last inequality follows, for instance, from $\e<3$ and $1/27>2/99$.
Therefore $F''(3)<0$ and $H''(3)>0$, excluding classical complete monotonicity of $F$ and the classical Bernstein property of $H$.

The same functions separate the ordinary definiteness conditions from the Jackson sign conditions.
For small $h>0$, choose $x_1=(3-h)/2$, $x_2=(3+h)/2$ and $c_1=1$, $c_2=-1$.
The quadratic form for the kernel $F(x+y)$ equals
\[
 F(3-h)-2F(3)+F(3+h)=h^2F''(3)+o(h^2)<0.
\]
Thus $F$ is not positive definite with respect to ordinary addition.
For the kernel $N_H$ in \eqref{eq:negativekernel}, the corresponding form is
\[
 -H(3-h)+2H(3)-H(3+h)=-h^2H''(3)+o(h^2)<0.
\]
Thus $H$ is not negative definite in that ordinary-addition sense.
\end{remark}

\begin{proposition}\label{prop:absolutecomposition}
Let $g$ be a classical completely monotone function on $(0,\infty)$, and let $f\in C^\infty(J)$ be absolutely monotone on an open interval $J$ containing the range of $g$; that is, $f^{(k)}\geq0$ on $J$ for every $k\geq0$.
Then $f\circ g$ is $q$-CM.
\end{proposition}
\begin{proof}
The classical composition property follows directly from the Fa\`a di Bruno formula: each term in $(f\circ g)^{(n)}$ has the form
\[
 C f^{(k)}(g(x))\prod_{j=1}^{k}g^{(m_j)}(x),
 \qquad C\geq0,\quad m_j\geq1,\quad \sum_{j=1}^{k}m_j=n.
\]
Multiplication by $(-1)^n$ makes every such term nonnegative.
Also $f\circ g\geq0$ by the order-zero assumption on $f$.
Thus $f\circ g$ is classically completely monotone, and Theorem~\ref{th:classicalcm} gives the conclusion.
\end{proof}

\begin{proposition}\label{prop:compositionlimits}
For a positive function $f$ on $(0,\infty)$, consider the conditions
\begin{enumerate}
\item[\textup{(i)}] $f$ is $q$-Bernstein;
\item[\textup{(ii)}] $g\circ f$ is $q$-CM for every $q$-CM function $g$ on $(0,\infty)$;
\item[\textup{(iii)}] $E_q(-kf)$ is $q$-CM for all sufficiently small $k>0$.
\end{enumerate}
Then \textup{(ii)} implies \textup{(i)}, and \textup{(iii)} implies \textup{(i)}.
Neither converse holds in general, even for positive smooth $f$.
Moreover, no strictly positive $f$ can satisfy the stronger requirement that $E_q(-kf)$ be $q$-CM for every $k>0$.
\end{proposition}
\begin{proof}
Assume \textup{(ii)} and choose $g(u)=\e^{-ku}$, $k>0$.
This function is $q$-CM by Theorem~\ref{th:classicalcm}.
For fixed $x$ and $n\geq1$, expansion at the finitely many points used by $\dq^n$ gives
\[
 0\leq(-1)^n\dq^n(\e^{-kf})(x)
 =k(-1)^{n-1}\dq^nf(x)+O_{x,n,f}(k^2).
\]
Dividing by $k$ and letting $k\downarrow0$ proves \textup{(i)}.
The implication \textup{(iii)}$\Rightarrow$\textup{(i)} is Proposition~\ref{prop:smallparameter}.

For failure of \textup{(i)}$\Rightarrow$\textup{(ii)}, use $q=1/10$ and the function $F$ in \eqref{eq:counterexamples}.
Take $f(x)=3+x$ and $g=F$.
Then $f$ is $q$-Bernstein and $g$ is $q$-CM, but \eqref{eq:multipleaverage}, applied at order two to $x\mapsto F(3+x)$, gives
\[
 \lim_{x\downarrow0}\dq^2(g\circ f)(x)
 =\dfrac{1+q}{2}F''(3)<0.
\]
Hence $g\circ f$ fails the second-order $q$-CM inequality near zero.

For failure of \textup{(i)}$\Rightarrow$\textup{(iii)}, let $f(x)=1+x$ and choose a fixed number $s$ satisfying
\begin{equation}\label{eq:negativerange}
 \dfrac{1}{1-q}<s<\dfrac{1}{q(1-q)}.
\end{equation}
In the product $E_q(-s)=\prod_{j\geq0}(1-(1-q)sq^j)$, the first factor is negative and all remaining factors are positive; their infinite product is positive.
Thus $E_q(-s)<0$.
For every $0<k<s$, the point $x_k=s/k-1$ is positive and
\[
 E_q(-kf(x_k))=E_q(-s)<0.
\]
Consequently $E_q(-kf)$ is not $q$-CM for any such $k$, although $f$ is $q$-Bernstein.

Finally, if $f(x_0)>0$ at any point, take $s$ as in \eqref{eq:negativerange} and $k=s/f(x_0)$.
Then $E_q(-kf(x_0))<0$.
This rules out the requirement for every $k>0$ for any strictly positive $f$.
\end{proof}

\begin{proposition}\label{prop:negativekernel}
Let $f$ be $q$-Bernstein on $\Xq$ and have a finite nonnegative limit at zero.
For a spectral representation \eqref{eq:qlevy}, supplied by Theorem~\ref{th:representation}, define
\begin{equation}\label{eq:formaltranslation}
 \mathcal T_f(x,y)=b+a(x+y)
 +\int_{\Sq}\bigl(1-K_t(x)K_t(y)\bigr)\,\mathrm d\mu(t),
 \qquad x,y\in\Xq.
\end{equation}
Then
\begin{equation}\label{eq:qnegativekernel}
 \mathcal N_f(x,y)=f(x)+f(y)-\mathcal T_f(x,y)
\end{equation}
is a positive definite kernel on $\Xq$.
Every classical Bernstein function, extended by its finite limit at zero, is negative definite with respect to ordinary addition.
The latter conclusion does not hold for all $q$-Bernstein functions.
\end{proposition}
\begin{proof}
The integral in \eqref{eq:formaltranslation} is finite, because $0\leq K_t\leq1$ and
\[
 0\leq1-K_t(x)K_t(y)
 \leq(1-K_t(x))+(1-K_t(y))
 \leq\min(1,tx)+\min(1,ty).
\]
Substituting \eqref{eq:qlevy} into \eqref{eq:qnegativekernel}, and cancelling the linear terms, yields
\[
 \mathcal N_f(x,y)
 =b+\int_{\Sq}(1-K_t(x))(1-K_t(y))\,\mathrm d\mu(t).
\]
For $x_1,\ldots,x_m\in\Xq$ and $c_1,\ldots,c_m\in\C$, its quadratic form is
\[
 b\left|\sum_{i=1}^m c_i\right|^2
 +\int_{\Sq}\left|\sum_{i=1}^m c_i(1-K_t(x_i))\right|^2\,\mathrm d\mu(t)
 \geq0.
\]
The integral is finite by Cauchy--Schwarz and the bounds
$(1-K_t(x_i))^2\leq1-K_t(x_i)\leq\min(1,tx_i)$.
This proves the asserted positivity.

For a classical Bernstein function, the ordinary negative-definiteness assertion is the standard consequence of \eqref{eq:classicallevy}; see \cite[Chapter~4]{Schilling}.
In that case the associated kernel is
\[
 N_f(x,y)=b+\int_{(0,\infty)}(1-\e^{-tx})(1-\e^{-ty})\,\mathrm d\mu(t).
\]
For the failure in the general $q$-Bernstein class, the function $H$ in Remark~\ref{rem:counterexample} provides a positive smooth counterexample.
\end{proof}

\medskip
\noindent\textit{Note on the scope of the revision.}
This article is a revised and corrected treatment of the paper by Ch. G. Kokologiannaki and V. Krasniqi \cite{KokologiannakiKrasniqi2014}.
It is important to distinguish \emph{correction} from \emph{generalization}.
In particular, certain statements from the 2014 article have been replaced by rigorously justified formulations with narrower domains or additional hypotheses, rather than extended in their original generality.
The spectral representation and its consequences in the present paper are established on the geometric lattice $\Xq$, with spectral measures supported on $\Sq$ and with the continuity and integrability assumptions explicitly imposed in the relevant results.
These restricted-domain conclusions must not be interpreted as corresponding assertions on arbitrary $q$-stable subsets of $(0,\infty)$ or on the whole positive half-line.
Any additional results stated here are to be understood only in the precise scope of their respective hypotheses.

The distinction is especially relevant to \emph{positive and negative definiteness}.
Proposition~\ref{prop:positivekernel} proves that the explicitly defined kernel $\mathcal K_g$ is positive semidefinite, and Proposition~\ref{prop:negativekernel} proves the same for the defect kernel $\mathcal N_f$, in both cases by an integral representation as a nonnegative mixture of rank-one kernels.
This does not imply that $g(x+y)$ is positive definite, or that a general $q$-Bernstein function is negative definite with respect to ordinary addition; the latter assertion is justified here only for classical Bernstein functions, and counterexamples show that the corresponding unrestricted $q$-claims fail.
Moreover, the formal $q$-addition identity in Remark~\ref{rem:addition} is an identity of formal powers and does not define an associative binary operation on $\Xq$.
The two-variable expression $\mathcal T_f$ is defined spectrally and is not claimed to arise from a point-valued semigroup translation.
Consequently, the positive semidefiniteness of $\mathcal K_g$ and $\mathcal N_f$ must not be taken to confer an unspecified semigroup structure, a convolution law, or semigroup-theoretic definiteness beyond the explicitly proved kernel statements.

\end{document}